\documentclass[11pt,reqno]{amsart}
\usepackage{amssymb,latexsym}
\usepackage{amsmath}
\usepackage{amsthm}

\numberwithin{equation}{section}
\newtheorem{theorem}{Theorem}[section]
\newtheorem{proposition}[theorem]{Proposition}
\newtheorem{lemma}[theorem]{Lemma}

\usepackage{multirow}

\theoremstyle{definition}

\newtheorem{remark}[theorem]{Remark}

\def\begeq{\begin{equation}}
\def\endeq{\end{equation}}

\allowdisplaybreaks

\begin{document}

\title{ New type of solutions for the Nonlinear Schr\"odinger  Equation in $\mathbb{R}^N$}
\author{Lipeng Duan}
\address{Lipeng Duan,
\newline\indent School of Mathematics and Statistics, Central China Normal University,
\newline\indent Wuhan 430079, P. R. China.
}
\email{ahudlp@sina.com}

\author{Monica Musso}
\address{Monica Musso,
\newline\indent Department of Mathematical Sciences, University of Bath,
\newline\indent Bath BA2 7AY, United Kingdom
}
\email{mm2683@bath.ac.uk}

\begin{abstract}
We construct a new family of entire solutions for the    nonlinear Schr\"odinger  equation
 \begin{align*}
 \begin{cases}
 -\Delta u+     V(y ) u = u^p, \quad u>0,   \quad \text{in}~ \mathbb{R}^N,
 \\[2mm]
u \in  H^1(\mathbb{R}^N),
\end{cases}
\end{align*}
where $p\in (1,      \frac{N+2}{N-2})$ and $N\geq 3$,   and $V (y)= V(|y|)$ is a positive bounded radial potential satisfying
$$
V(|y|) =  V_0 + \frac{a}{|y|^m} + O( \frac{1}{|y|^{m+\sigma}} ),     \quad {\mbox {as}} \quad |y| \to \infty ,
$$
for some fixed constants $V_0, a, \sigma >0$, and $m>1$. Our solutions are different from the ones obtained in \cite{wei&yan} and have
strong analogies with the doubling construction of entire finite energy sign-changing solution for the Yamabe equation in \cite{MaMo}.

\vspace{4mm}


\vspace{2mm}

\end{abstract}

\date{\today}


\maketitle

 \section{introduction}

This paper is devoted to the construction of  solutions to the  following    nonlinear elliptic problem
 \begin{align}
 \begin{cases}\label{original}
 -\Delta u+     V(y ) u = u^p,  \quad u>0,  \quad \text{in}~ \mathbb{R}^N,
 \\[2mm]
u \in  H^1(\mathbb{R}^N),
\end{cases}
\end{align}
where $ 1< p< \frac{N+2}{N-2}$, $N\geq 3$, and $V$ is a bounded radially symmetric potential with $V(y) \geq V_0 >0$.
This problem arises when looking for standing waves solutions
\begin{align*}
\psi(y,t) = e^{{  \bf{ i } } \lambda t } u(y).
\end{align*}
to the   time-dependent  Schr\"odinger equations
  \begin{align}\label{original1}
-{  \bf{ i } } \frac{\partial \psi}{\partial t}  = \Delta \psi - V_1(y) \psi + |\psi|^{p-1} \psi,     \quad\text{for} ~ (y,t) \in \mathbb{R}^N\times \mathbb{R},
 \end{align}
where $ {  \bf{ i } }$ is the imaginary unit and   $   V(y) = V_1(y) +\lambda$.  This kind of equation arises in many applications, for instance in  nonlinear optics, plasma physics,  quantum mechanics, or condensed matter physics.

 \vspace{1mm}
The energy functional  associated to  \eqref{original} is given by
\begin{equation}\label{I1}
I(u) = \frac 12 \int_{  \mathbb{R}^N}   \Big\{ |  \nabla u |^2 + V(y) u^2\Big\}  - \frac{1}{p+1} \int_{\mathbb{R}^N}  |u|^{p+1}.
\end{equation}
This functional posses infinitely many critical points in the class of radially symmetric functions, as a direct consequence of the
Ljusternik-Schnirelmann theory. Nonetheless, it is not clear that these critical points are solutions to \eqref{original} since they may not need to be positive. On the other hand, using
 the  concentration compactness  theorem  one can prove that if the potential $V=V(y)$ satisfies $ V(y ) \leq V_0 $ for all $y\in  \mathbb{R}^N$, then Problem  \eqref{original} has a least energy solution;  while if $V(|y|) \geq V_0$ and $V(|y|) \neq V_0$,  then \eqref{original} does not have a least energy solution, \cite{Lions1,Lions2}, and solutions have high energy.
In \cite{CDS},    the authors  showed the existence of infinitely many sign-changing solutions for Problem  \eqref{original}  under the assumptions that
$\varliminf_{|y|  \rightarrow + \infty      } V(y) = V_0 >0$  with a suitably rate.  They used sophisticated variational arguments which allow them to treat the problem without    requiring  any  symmetry property on  $ V(y)$. We refer the reader to the articles  \cite{AMN0},  \cite{AMN1}, \cite{DWK1}, \cite{FW},  \cite{LinNiWei}, \cite{WY1} , \cite{WXY}     and the references therein for related results, also for the associated
 singularly perturbed problem.

 \vspace{1mm}

In this paper we assume that the potential $V$ is radial and
 satisfies  the following decay condition at infinity:  there exist $a, \sigma >0$, and $m>1$ so that
\begin{align}\label{H1}
V(|y|) =  V_0 + \frac{a}{|y|^m} + O_k( \frac{1}{|y|^{m+\sigma}} ),     \quad {\mbox {as}} \quad |y| \to \infty .\tag{H1}
\end{align}
Observe that a direct scaling argument enables us to just consider
 $V_0=1$. Under these same assumptions, Wei and Yan in \cite{wei&yan} used a constructive method to produce infinitely many non-radial solutions to \eqref{original} with high energy. Roughly speaking, these solutions are obtained by gluing together a large number of copies of building blocks given by the  {\it bump} $U$ which is the only positive radially-symmetric solution to
\begin{align} \label{Schr}
 \begin{cases}
 -\Delta u + u= u^p,       \qquad  &\text{in}~ \mathbb{R}^N,
 \\[2mm]
u(y)    \rightarrow 0,      \quad &\text{as}~ |y| \rightarrow \infty.
\end{cases}
\end{align}
 More precisely, for any large integer $k$ the authors constructed a solution $u_k$ looking like   a sum of a $k$   {\it bumps}  $U(x-x_j^*)$,
 \begin{equation}\label{uuk}
 u^*_k (y) \sim \sum_{j=1}^k U(y-x_j^* )
 \end{equation}
 where  the  location points
 $x_j^*$  are distributed along the vertices of a regular $k$-polygon
$$ x_j^* = r \Big(      \cos{ \frac{2(j-1) \pi }{k}     },      \sin{ \frac{2( j-1) \pi }{k}     },  0,   \cdots, 0 \Big), \quad\text{for}~ j= 1,\cdots, k$$
with large radius $r \sim k \log k$ as $k \to \infty$. These solutions respect the polygonal symmetry in the $(x_1,x_2)$-plane, and are radially symmetric in the other variables. This construction is stable in the sense that solutions with similar profile exist for other potential functions close to $V$ but not necessarily radially symmetric \cite{DWY}. In fact, the solutions constructed in \cite{wei&yan} happen to be {\it non-degenerate} in the class of functions sharing the same symmetry, in the sense that the linearized operators around these solutions are invertible. This fact has been recently proven in \cite{GMPY} with the use of  a refined version of the Pohozaev identity in the spirit of \cite{GMPY1}. This property makes it possible to use the functions $u_k^*$ as new building blocks   to  generate new constructions, see \cite{GMPY}.

  \vspace{1mm}
A natural question to ask is whether there are other type of {\it building blocks} for \eqref{original}, which are not equal to the ones constructed in \cite{wei&yan} nor cannot be obtained by gluing a number of them together. The purpose of this paper is to give an  answer of this question with a new family of solutions to \eqref{original} which have a more complex  concentration structure.

Let $N\geq 3$,  $k$ be an integer and introduce the points
\begin{align} \label{overunder}
 \begin{cases}  \overline{x}_j= r \Big(    \sqrt {1-h^2 }  \cos{ \frac{2(j-1) \pi }{k}     },    \sqrt {1-h^2 } \sin{ \frac{2( j-1) \pi }{k}     }, h,    {\bf{  0  } }  \Big), \quad j= 1, \cdots, k,
 \\[5mm]
  \underline{x}_j= r \Big(    \sqrt {1-h^2 }    \cos{ \frac{2(j-1) \pi }{k}     },   \sqrt {1-h^2 }  \sin{ \frac{2( j-1) \pi }{k}     }, - h,       {\bf{  0  } }   \Big), \quad j= 1, \cdots, k,
 \end{cases}
\end{align}
where $ {\bf{  0  } }$ is the zero vector in $\mathbb{R}^{N-3}$. The parameters $h$ and $r$ are positive numbers and are chosen in the range
 \begin{equation}\label{par1}
 h \in \big[ \alpha_0\frac{1}{k}, \alpha_1\frac{1}{k}\big]   , \quad  r\in [\beta_0 k \ln k, \beta_1 k \ln k]
 \end{equation}
 for $\alpha_0, \alpha_1, \beta_0, \beta_1$ fixed positive constants, independent of $k$.
We  define the approximate solution  as
\begin{align}\label{Wrh1}
W_{r,h} (y) = \sum_{j=1}^kU_{\overline{x}_j } (y) + \sum_{j=1}^k U_{\underline{x}_j } (y) ,
\end{align}
where $U_{\overline{x}_j } (y) = U (y-\overline{x}_j )  $ and $ U_{\underline{x}_j } (y) = U(y-\underline{x}_j )$ and $k$ would
be sufficiently large.        In this paper, we will prove that  from any $k$ large enough problem \eqref{original}  has a new type of solutions $u_k $ which have the form
\begin{equation}\label{uk} u_k (y) \sim W_{r,h} (y) ,
\end{equation} as $k \to \infty$. Our solutions will have  polygonal symmetry in the $(x_1,x_2)$-plane, will be even in the $x_3$ direction and  radially symmetric in the variables $x_4, \ldots , x_N$. In particular, they do not belong to the same class of symmetry as the solutions built in \cite{wei&yan}. Our solutions are thus different from the ones obtained in \cite{wei&yan} and have
strong analogies with the doubling construction of entire finite energy sign-changing solution for the Yamabe equation in \cite{MaMo}.

  \vspace{1mm}
Observe that, if we take $h=0$ in \eqref{overunder}, then $\overline{x}_j= \underline{x}_j$ for any $j$ and the two constructions  \eqref{uuk} and \eqref{uk}-\eqref{Wrh1} are the same.  Under our assumptions on the range for the parameters $h$ and $r$ in \eqref{par1}, the two constructions \eqref{uuk} and \eqref{uk}-\eqref{Wrh1} are different.     A main  distinction between    the present  construction and the one in  \cite{wei&yan}    is that  there are two parameters $r, h$ to choose in the locations $ {\underline{x}_j },    {\overline{x}_j }   $ of the  bumps in \eqref{Wrh1}.

We will discuss in details our result in the following subsection.

Throughout this paper,  we employ  $C,    C_j,   \text{or}~ \sigma, \sigma_j,  \tau, \tau_j,   j=0,1,2,\cdots$ to denote certain constants.
Furthermore,  we   also     employ the common  notation    by writing $O_k(f(r,h)), o_k(f(r,h))  $ for the functions which satisfy
\begin{align*}
\text{if} \quad  g(r,h) \in O_k(f(r,h)) \quad \text{then }\quad  {\lim_{k \to +\infty}} \Big|  \,  \frac{g(r,h)}{f(r,h)} \, \Big|  \leq C< + \infty,
\end{align*}
and
\begin{align*}
\text{if} \quad  g(r,h) \in o_k(f(r,h)) \quad \text{then }\quad  {\lim_{k \to +\infty}}  \frac{g(r,h)}{f(r,h)} =0,
\end{align*}
  in the present paper.

\vspace{3mm}
 \subsection{Main result and scheme of the proof }

\textsc{}
\vspace{2mm}

For     $ j= 1,\cdots, k$, we divide    $\mathbb{R}^N$ into $k$ parts:
\begin{align*}
\Omega_j : = &\Big\{  y = (y_1, y_2, y_3, y'') \in \mathbb{R}^3 \times \mathbb{R}^{N-3}\nonumber
 \\[1mm]
& \qquad
:       \langle  \frac{ (y_1, y_2)}{|(y_1, y_2) |},           ( \cos{ \frac{2(j-1) \pi }{k}     },    \sin{ \frac{2( j-1) \pi }{k}     } )   \rangle_{\mathbb{R}^2 }\geq \cos{ \frac \pi k}          \Big\}.
\end{align*}
where $ \langle   ,  \rangle_{\mathbb{R}^2 } $  denote  the dot product in $\mathbb{R}^2$.
 For $\Omega_j$,   we divide  it into  two parts:
\begin{align*}
\Omega_j^+ = & \Big\{  y:    y = (y_1, y_2, y_3, y'')   \in   \Omega_j, y_3\geq0 \Big\},
\end{align*}
\begin{align*}
\Omega_j^- = & \Big\{  y:  y = (y_1, y_2, y_3, y'')   \in   \Omega_j, y_3<0 \Big\}.
\end{align*}
We see that
\begin{align*}
\mathbb{R}^N =  \cup_{j=1}^k   \Omega_j,   \quad  \Omega_j =     \Omega_j^+    \cup    \Omega_j^-
\end{align*}
and  the  interior of
\begin{align*}
 \Omega_j  \cap   \Omega_i,        \quad   \Omega_j^+    \cap    \Omega_j^-
\end{align*}
   are   empty sets for $i\neq j$.

\medskip

We now define the  symmetric Sobolev space:
\begin{align*}
H_s = \Bigg\{  & u : u \in H^1(\mathbb{R}^N),    \text{ $u$   is even in   $ y_\ell, \ell = 2, 4, 5, \cdots, N,$}  \quad   u \big( \sqrt{y_1^2+y_2^2} \cos \theta, \sqrt{y_1^2+y_2^2} \sin \theta,   \nonumber
\\
       &  y_3, y''       \big)    =  u \big( \sqrt{y_1^2+y_2^2} \cos { \big( \theta+   \frac{2j\pi}{k}  \big)  }, \sqrt{y_1^2+y_2^2} \sin { \big( \theta+   \frac{2j\pi}{k}  \big)  }, y_3, y''       \big)         \Bigg\}.
\end{align*}
where $ \theta =  \arctan{\frac {y_2}{ y_1} }$.

\vspace{2mm}
   In this paper, we always assume
\begin{align}\label{H2}
 (r,h) \in \mathbb{S}_k   & =: \big[  \big(\frac{m}{2 \pi    }- \beta  \big)   k \ln k,   \,       \big(\frac{m}{2 \pi    }+ \beta  \big)   k \ln k\big]  \nonumber
 \\[1mm]
 & \qquad
 \times \big[   \big( \frac{ \pi \, (m+2) }{  m   } - \alpha \big) \frac{1}{k},    \,  \big( \frac{ \pi \, (m+2) }{  m    }  + \alpha \big) \frac{1}{k} \big],
\end{align}
for some $ \alpha, \beta>0$  small, and independent of $k$.  We refer to Remark \ref{remark6} for a discussion on the assumption \eqref{H2}  for $(r,h)$.


Our main result is the following

\vspace{3mm}
\begin{theorem}\label{main1}
Suppose that  $ V(|y|) $ satisfies  \eqref{H1} and the parameters $(r,h)$ satisfies   \eqref{H2}.  Then there is an integer $k_0 $,   such that for all integer $k\geq k_0$,
\eqref{original} has a solution $u_k $ of the form
\begin{align}\label{u_k}
u_k = W_{r_k, h_k}(y) + \omega_k(y),
\end{align}
where $ \omega_k \in    H_s,  (r_k,s_k) \in \mathbb{S}_k  $ and  $ \omega_k $ satisfies
\begin{align*}
  \int_{\mathbb{R}^N }  \big( | \nabla {\omega}_k |^2 + V(y) |\omega_k|^2 \big)  \to 0, \quad {\mbox {as}} \quad k \to \infty.
\end{align*}
\qed
\end{theorem}

\medskip
We will prove  Theorem \eqref{main1} by using the Lyapunov- Schmidt reduction technique adapted to our context as developed in \cite{wei&yan}. Let  us briefly sketch the scheme of the proof. A critical point $u$ of the energy functional $I$  defined in \eqref{I1} corresponds to a solution for \eqref{original}.   Our solution $u$ will have the form    $u= W_{r,h}   + \phi$.

As we know,    \eqref{Schr}  has a unique  positive solution $U$, which is radially symmetric and
$$  \lim_{|y| \rightarrow +\infty}  U(y)  e^{|y|} y^{\frac{N-1}{2}} = C<+\infty,     \quad\text{and } ~ \lim_{|y| \rightarrow +\infty} \frac{U(y)}{U'(y)}=-1.    $$
Moreover,  $U$ is non-degenerate, in the sense that the kernel of the linear operator  $-\Delta + \mathbb{I} -
pU^{p -1} $  in $H^1(\mathbb{R}^N) $ is spanned by  $ \{ \frac{\partial U}{\partial {y_1}  }, \cdots, \frac{\partial U}{\partial {y_N} } \}$.   For more details about  \eqref{Schr}, reader can refer to  \cite{Kwong}, \cite{NWM1}.

We define
\begin{align*}
J(\phi ) =     I (W_{r,h}   + \phi ),       \quad  \phi \in    \mathbb{E}.
\end{align*}
The space $\mathbb{E}$ is given as follows. \medskip
 For $ j = 1, \cdots, k, $    we define
 \begin{align*}
 \overline{ \mathbb{Z}}_{1j}  =    \frac{\partial U_ {\overline{x}_j} }{\partial r}, \quad \underline{ \mathbb{Z}}_{1j}  =    \frac{\partial U_ {\underline{x}_j} }{\partial r}, \quad  \overline{ \mathbb{Z}}_{2j}  =    \frac{\partial U_ {\overline{x}_j} }{\partial h},  \quad \underline{ \mathbb{Z}}_{2j}  =    \frac{\partial U_ {\underline{x}_j} }{\partial h}.
 \end{align*}
We  define the constrained  space
\begin{align}\label{SpaceE}
\mathbb{E}  = \Big\{ v:   v\in H_s,      \quad \int_{{\mathbb{R}}^N } & U_{\overline{x}_j}^{p-1}   \overline{ \mathbb{Z}}_{\ell j} v = 0  \quad\text{and}    \nonumber
\\[2mm]
  \quad  & \int_{{\mathbb{R}}^N }  U_{\underline{x}_j}^{p-1}   \underline{ \mathbb{Z}}_{\ell j}  v = 0, \quad j = 1, \cdots, k, \quad \ell = 1,2      \Big\}. \end{align}
  The space $ \mathbb{E}  $   will be endowed with the norm
$$ \| v\|^2   =   \langle v, v\rangle,         \quad ~ v \in \mathbb{E},   $$
where
$$ \langle v_1, v_2\rangle  =   \int_{{\mathbb{R}}^N } \big (     \nabla {v_1} \nabla {v_2}  + V (|y|) v_1 v_2  \big) \quad v_1, v_2 \in        \mathbb{E}.$$
Expand $ J(\phi )$ as follows:
\begin{align*}
J(\phi ) = J(0) + {  \bf{ l } }  (\phi) + \frac 12  \langle{\bf{L}} \phi,  \phi\rangle  -  {  \bf{ R } }(\phi),    \quad \phi \in   \mathbb{E},
\end{align*}
where
\begin{align*}
{  \bf{ l } }  (\phi)  = \sum_{i =1}^k   \int_{  \mathbb{R}^N} \big( V(|y|) - 1      \big)  \big( U_{ \overline{x}_i  } + U_{\underline{x}_i  }\big)  \phi +   \int_{  \mathbb{R}^N}  \Big(   \sum_{i =1}^k\big( U_{ \overline{x}_i  } + U_{\underline{x}_i  }\big)^p -  W_{r,h}^p    \Big) \phi,
\end{align*}
and
\begin{align*}
{  \bf{ R } }(\phi) = \frac{1}{p+1}  \int_{  \mathbb{R}^N} \Big(| W_{r,h}+\phi|^{p+1} -  | W_{r,h} |^{p+1} - (p+1)  W_{r,h}^{p}-   \frac 12 (p+1) pW_{r,h}^{p-1} \phi^2  \Big).
\end{align*}
Furthermore ${\bf{L}} $ is  a   linear operator  from $  \mathbb{E} $ to $ \mathbb{E}$,  which satisfies
$$     \langle{\bf{L}} v_1,  v_2 \rangle  =     \int_{{\mathbb{R}}^N } \big (     \nabla {v_1} \nabla {v_2}  + V (|y|) v_1 v_2  -      p W_{r,h}^{p-1} v_1v_2  \big),   \quad    \text{for all }   ~ v_1, v_2 \in       \mathbb{E}.$$
Since  $  W_{r,h}^{p-1}$ is bounded and has the symmetries of the space $H_s$,     we   can easily
 prove  $  {\bf L}  $   is a  bounded   linear operator from $  \mathbb{E}$ to $\mathbb{E}$.
We will show that $ {  \bf{ l } }  (\phi)$ is a bounded linear functional in $ \mathbb{E}$.  Thus  there is an ${  \bf{ l } }_k\in  \mathbb{E}$, such that  $$   {  \bf{ l } }  (\phi)=    \langle  {  \bf{ l } }_k, \phi  \rangle.$$
Then  a critical point of $ J(\phi )$  is also  a solution of
\begin{align}\label{equivalentequat}
{  \bf{ l } }_k +{\bf{L}} \phi+  {  \bf{ R } }' (\phi) = 0.
\end{align}
Thus a function $u$ of the form $u= W_{r,h}   + \phi$ will be a solution to \eqref{original} if   $ \phi$ is a solution of \eqref{equivalentequat}.
The strategy now consists in first showing that, for any $ (r,h) \in \mathbb{S}_k $, there exists a function  $\phi_{r,h} $  solution  of  \eqref{equivalentequat}  in the  space $\mathbb{E}$ (see Proposition \ref{propos2}).   A second step will be to reduce the problem of finding a critical point $u$ of  $I(u)$ to  the problem of finding a stable critical point $(r^*,h^*)$ of the  function $$    F(r,h)=  I (W_{r,h}   + \phi_{r,h}  ).$$
We will show that such a critical point exists in the set $\mathbb{S}_k$. The qualitative property of the solutions  follows by their construction.


%

 \vspace{3mm}

\noindent{\textbf{Plan of the paper}}\\
We organize the paper as follows.      In section 2, we will give the proof of  Theorem \ref{main1}.    In the Appendices, we will give some useful Lemmas,  Propositions and  the  details  for the energy of approximate solution  expansion.

\section{   Proof of Theorem \ref{main1}}

We  first give    expansion  for the energy of  approximate solution.
\medskip
\begin{proposition}\label{func}
 For all $(r,h) \in \mathbb{S}_k $, there exist some  small constant $\sigma>0$  such that
\begin{align}
 I(W_{r,h} )
   =   k \Big(      \frac{A_1}{r^m} + A_2    &   - 2B_1     e^{-     2 \pi \sqrt{1-h^2 }    \frac{r}{k}    } -B_1 e^{-2    rh}   + O_k(\frac{1}{r^{m+\sigma}}) \Big)      \nonumber
   \\[1mm]
&  \quad     + k O_k(e^{- 2(1+\sigma)    rh}  )+ k O_k(     e^{-     2(1+\sigma )\sqrt{1-h^2 }    \frac{r}{k}    }     ),
\end{align}
where
\begin{align}\label{A1A2}
A_1 = a \int_{  \mathbb{R}^N}  U^2,   \quad     A_2 =   \big(1- \frac{2}{p+1} \big) \int_{   \mathbb{R}^N}     U^{p+1},  \quad     B_1=     \int_{  \mathbb{R}^N}  U^p e^{-y_1}. \end{align}
\end{proposition}

{\textit{Proof of Proposition \ref{func}  } }:
The proof  of Proposition \ref{func}  is  delayed to   Appendices. \qed

\medskip
The next lemma gives the existence and boundness of inverse operator of ${\bf{L}}$ in    $ \mathbb{E}. $

\medskip
\begin{lemma}\label{lemmainverse}
There is a constant $ \rho > 0$, independent of $k$, such that for any $ (r,h)\in  \mathbb{S}_k$
$$  \|{\bf{L}} v\| \geq \rho \| v\|, \quad v\in \mathbb{E}. $$
\end{lemma}
{\textit{Proof of Lemma \ref{lemmainverse}}}:
We prove by contradiction.   Suppose that  when  $k\rightarrow +\infty$
there exist  $h_k, r_ k  \in   \mathbb{S}_k,       v_k \in  \mathbb{E}  $ satisfying
\begin{align*}
\|{\bf{L}} v_k\| =  o_k(1) \| v_k\|.
\end{align*}
Then easily
\begin{align*}
\langle{\bf{L}} v_k,      \varphi  \rangle     =     o_k(1) \| v_k\| \,  \|    \varphi    \|,    \quad \forall\,  \varphi \in  \mathbb{E}.
\end{align*}
Similar to \cite{wei&yan}, we assume $\| v_k\|^2   =k$.     Using the symmetric property   of $ v_k, \varphi$, we can get
\begin{align} \label{operatorL}
\langle{\bf{L}} v_k,      \varphi  \rangle   &=   \int_{{\mathbb{R}}^N } \big (     \nabla {v_k } \nabla {\varphi }  + V (|y|) v_k  \varphi  -      p W_{r,h}^{p-1} v_k\varphi  \big) \nonumber
 = k     \int_{ \Omega_1}  \big (     \nabla {v_k } \nabla {\varphi }  + V (|y|) v_k  \varphi  -      p W_{r,h}^{p-1} v_k\varphi  \big) \nonumber
\\[1mm]
& =o_k(1) \| v_k\| \,  \|    \varphi    \|   =o(\sqrt k)  \,  \|    \varphi    \|        ,
\end{align}
and
\begin{align}\label{norm}
  \int_{ \Omega_1}  \big (     | \nabla {v_k } |^2   + V (|y|) v_k^2    = 1.
\end{align}
Inserting  $  \varphi  =  v_k $ into \eqref{operatorL},    we can obtain immediately
\begin{align*}
  \int_{ \Omega_1}  \big (   |  \nabla {v_k } |^2 + V (|y|) v_k^2 -      p W_{r,h}^{p-1} v_k^2  \big)    = o_k(1).
\end{align*}

We denote
$$ \overline{v}_k (y) = v_k (y+  \overline{x}_1 ).   $$
For the sequence $ \overline{v}_k (y) $, we can prove that  $ \overline{v}_k (y) $ is bounded in $  H^1_{loc}(\mathbb{R}^N). $
In fact,
for  any $ R>0$ ,      since $ | \overline{x}_2 -  \overline{x}_1| = 2r  \sqrt{1-h^2} \, \sin{      \frac \pi k } \geq  \frac{m}{4} \ln k,  $
we can choose $k$ large enough such that  $ B_R(   \overline{x}_1) \subset \Omega_1$.       As a result, we have
\begin{align}
  \int_{  B_R(   0 ) }  \big ( \,     | \nabla { \overline{v}_k } |^2   + V (|y  |)  \overline{v}_k^2 \big ) & =  \int_{  B_R(   \overline{x}_1)}  \big (  \,    | \nabla {  v_k } |^2   + V (|y- \overline{x}_1  |) { v}_k^2\, \big )  \nonumber
  \\[2mm]
    & \leq     \int_{ \Omega_1}  \big (\,      | \nabla {{ v}_k } |^2   + V (|y- \overline{x}_1|)  { v}_k^2 \,  \big )\leq 1.
\end{align}
 So we can conclude
 \begin{align}\label{condition0}
 \overline{v}_k \rightarrow   \bar{ v }\quad  \text{ weakly in } H^1_{loc}(\mathbb{R}^N),
 \end{align}
 and  \begin{align}\label{condition5}
 \overline{v}_k \rightarrow  \overline{v} \quad  \text{ strongly in } L^2_{loc}(\mathbb{R}^N).
 \end{align}
Since $\overline{v}_k$ is even in $y_d, d = 2, 4,  \cdots, N$, then $ {\bar v}$ is even in  $y_d$.

 From the orthogonal conditions for functions of $\mathbb{E} $
\begin{align*}
\int_{{\mathbb{R}}^N }  U_{\overline{x}_1}^{p-1}  \frac{\partial U_ {\overline{x}_1} }{\partial r } \,  {v}_k  = 0
\end{align*}
and the identity
\begin{align*}
\frac{\partial  U_ {\overline{x}_1} }     {  \partial r } =   \sqrt{1-h^2} \,     \frac{ \partial U_ {\overline{x}_1} }{\partial {y_1} }     + h \,   \frac{ \partial U_ {\overline{x}_1} }{\partial {y_3} },
\end{align*}
we can get
\begin{align} \label{condition1}
  \sqrt{1-h^2} \,     \int_{{\mathbb{R}}^N }  U_ {\overline{x}_1}^{p-1}       \frac{ \partial U_{\overline{x}_1} }{\partial {y_1} }     {v}_k       +  h \,    \int_{{\mathbb{R}}^N }  U_ {\overline{x}_1}^{p-1}       \frac{ \partial U_{\overline{x}_1} }{\partial {y_3} }    {v}_k =0.
  \end{align}
%
Similarly, combining
\begin{align*}
\int_{{\mathbb{R}}^N }  U_{\overline{x}_1}^{p-1}  \frac{\partial U_ {\overline{x}_1} }{\partial h } \,   {v}_k  = 0,
\end{align*}
and
\begin{align*}
\frac{\partial U_ {\overline{x}_1} }{\partial h } =  -  \frac{ h r}{ \sqrt{1-h^2} }  \frac{ \partial U_ {\overline{x}_1} }{\partial {y_1} }+ r \frac{ \partial U_ {\overline{x}_1} }{\partial {y_3} },
\end{align*}
we can get
\begin{align} \label{condition2}
  \frac{ h }{ \sqrt{1-h^2} }    \int_{{\mathbb{R}}^N }  U_ {\overline{x}_1}^{p-1}       \frac{ \partial U_{\overline{x}_1} }{\partial {y_1} }     {v}_k    - \int_{{\mathbb{R}}^N }  U_ {\overline{x}_1}^{p-1}       \frac{ \partial U_{\overline{x}_1} }{\partial {y_3} }     {v}_k =0.
\end{align}
From \eqref{condition1}, \eqref{condition2},  we have
\begin{align*}
\int_{{\mathbb{R}}^N }  U^{p-1}  \frac{\partial U  }{\partial y_1 } \,  \overline{v}_k   = \int_{{\mathbb{R}}^N }  U^{p-1}  \frac{\partial U  }{\partial y_3 } \,  \overline{v}_k =0.
\end{align*}
 Letting  $k \rightarrow   + \infty$,  we obtain
 \begin{align}\label{condition3}
 \int_{{\mathbb{R}}^N }  U^{p-1}       \frac{ \partial U }{\partial {y_1} }     \overline{v}    = \int_{{\mathbb{R}}^N }  U^{p-1}       \frac{ \partial U }{\partial {y_3} }     \overline{v}    =0.
 \end{align}


Next, we will show that  $\overline{v}$ is a solution of
\begin{align}\label{eqs0}
- \Delta \phi   + \phi   - p U^{p-1} \phi   =0, \quad \text{in }~ \mathbb{R}^N.
\end{align}
 We define  the constrained  space as:
  \begin{align*}
  {\tilde E}^+ = \Bigg\{     \phi : \phi \in H^1( \mathbb{R}^N ),       \int_{{\mathbb{R}}^N}   U^{p-1}       \frac{ \partial U }{\partial {y_1} }     \phi= \int_{{\mathbb{R}}^N }  U^{p-1}       \frac{ \partial U }{\partial {y_3} }    \phi =0
    \Bigg\}.
  \end{align*}
  For the proof of \eqref{eqs0},
  we first give a  claim.\\
  \textbf{Claim 1}:  $\overline{v}  $ is a solution of
  \begin{align*}
- \Delta \phi   +\phi   - p U^{p-1} \phi  =0, \quad \text{in }~   {\tilde E}^+.
\end{align*}

 Now we give the proof  of the   \textbf{Claim 1}.

For any $R>0 $, let  $\phi \in C_0^\infty \big (B_R(0) \big) \cap   {\tilde E}^+  $ which is even in $y_d, d = 2, 4, \cdots, N$.  Then denote $$ \phi_k (y)=: \phi (y -  \overline{x}_1  )  \in  C_0^\infty \big (B_R(  \overline{x}_1) \big).$$
 Inserting  $ \phi_k (y) =\varphi  $ into \eqref{operatorL} and   combining \eqref{condition0}, \eqref{condition5} and Lemma \eqref{lemma1},  we can get
 \begin{align} \label{eqs1}
 \int_{\mathbb{R}^N }  \Big( \nabla  \overline{v} \nabla\phi +  \overline{v} \phi - p U^{p-1}   \overline{v}   \phi \Big) = 0.
 \end{align}
 On the other hand,  since $ \overline{v}$ is even in $y_d, d = 2, 4, \cdots, N$,    then by using symmetric conditions  we can conclude   that \eqref{eqs1}   is valid for all functions  $\phi  \in C_0^\infty \big ( B_R(  \overline{x}_1 \big) \cap   {\tilde E}^+  $ which is odd in $y_d, d = 2, 4, \cdots, N$.    Hence \eqref{eqs1} is hold for all functions $\phi  \in C_0^\infty \big ( B_R(  \overline{x}_1 \big) \cap   {\tilde E}^+  $.
Density argument implies that,
\begin{align} \label{eqs2}
 \int_{\mathbb{R}^N }  \Big( \nabla  \overline{v} \nabla\phi +  \overline{v} \phi - p U^{p-1}    \overline{v}  \phi \Big) = 0,   \quad \text{for all} ~ \phi \in {\tilde E}^+.
 \end{align} The proof of \textbf{Claim 1} is completed.

  Combining  \textbf{Claim 1} and the fact that  \eqref{eqs0} is hold  for   $ \phi = \frac{\partial U}{\partial y_1} $ and $ \phi = \frac{\partial U}{\partial y_3} $,  then we get
  \begin{align}
 \int_{\mathbb{R}^N }  \Big( \nabla  \overline{v} \nabla\phi +  \overline{v} \phi - p U^{p-1}    \overline{v}  \phi \Big) = 0,   \quad \text{for all} ~ \phi \in  H^1(\mathbb{R}^N),
  \end{align}
which is \eqref{eqs0}. By using the Non-degeneracy results for $U$ and combining $ \overline{v} $ is even in   $y_d, d = 2, 4, \cdots, N$,     we have
  \begin{align}  \label{condition4}
  \overline{v} = c_1  \frac{\partial U}{\partial y_1}+  c_2 \frac{\partial U}{\partial y_3},
\end{align}
 for some universal constants $ c_1, c_2$.   Combining \eqref{condition3},\eqref{condition4}, we have
  \begin{align*}
c_1= c_2 =0.
\end{align*}
Thus we have
\begin{align}\label{v0} \overline{v} =0.
\end{align}
The direct result of \eqref{condition5} and \eqref{v0}  is that
\begin{align*}
  \int_{ B_R(  \overline{x}_1 ) }     v_k^2  =  o_k(1).
\end{align*}

From  Lemma \eqref{lemma1}, we have
$
W_{r,h} \leq C  e^{-(1-\eta) | y -   \overline{x}_1  | }.
$
Then we can get, taking $R$ large enough,
\begin{align*}
o_k(1) & =   \int_{ \Omega_1}  \big (   |  \nabla {v_k } |^2 + V (|y|) v_k^2 -      p W_{r,h}^{p-1} v_k^2  \big)  \nonumber
\\[1mm] & =  \int_{ \Omega_1}  \big (   |  \nabla {v_k } |^2 + V (|y|) v_k^2    - \int_{ \Omega_1\setminus  B_R(  \overline{x}_1 ) }    p W_{r,h}^{p-1} v_k^2    -  \int_{ B_R(  \overline{x}_1 ) }   p W_{r,h}^{p-1} v_k^2  \nonumber
\\[1mm] & =  \int_{ \Omega_1}  \big (   |  \nabla {v_k } |^2 + V (|y|) v_k^2 + o_k(1)+  O_k(e^{-(1-\eta)R  })  \int_{ \Omega_1}  |   {v_k } |^2 \nonumber
\\[1mm]
&\geq  \frac 12 \int_{ \Omega_1}  \big (   |  \nabla {v_k } |^2 + V (|y|) v_k^2 + o_k(1),
\end{align*}
which is a contradiction to \eqref{norm}.  We complete the proof of Lemma \eqref{lemmainverse}.
\qed


%
%

\medskip
Now we give the following Proposition  which is crucial in the sequel.

\medskip
\begin{proposition}\label{propos2}
For  any $k$ sufficiently large enough, there exist  a $C^1  $ map  $\Phi: \mathbb{S}_k   \rightarrow  \mathbb{E} $  such that   $ \Phi(r,h)  = \phi_{r,h} (y) \in    \mathbb{E}$, and
\begin{align}\label{critiacl}
J'(\phi_{r,h} ) = 0, \quad  \text{in } ~ \mathbb{E}.
\end{align}
Moreover, there is a small $\sigma>0$, such that
\begin{align*}
\|\phi_{r,h} \| \leq  \frac{C}{k^{  \frac{m-1}{2} + \sigma } }.
\end{align*}
\end{proposition}
{\textit{Proof of Proposition \ref{propos2}     }}:   The proof follows from a standard   technique,    which  is based on the contraction mapping theorem.  For the proof,   it's a     slight modification  of the proof of Proposition 2.2 in \cite{wei&yan}.  Here we omit for concise.
\qed

\medskip
Next we will give the estimate for $   {  \bf{ l } }_k$. The estimate will play a role in carrying out the reduction argument in  the proof of  Theorem \ref{main1}.
\vspace{1mm}

\begin{lemma}\label{lemmaestimate}
 For all $(r,h) \in \mathbb{S}_k $, there is a small $\sigma>0$, such that
\begin{align}\label{estimatelk}
\|  {  \bf{ l } }_k\| \leq    \frac{C}{k^{  \frac{m-1}{2} + \sigma } }.
\end{align}
\end{lemma}
{\textit{Proof of Lemma \ref{lemmaestimate}}}:    Recall  that for any $\phi \in \mathbb{E}$
\begin{align} \label{definlk}
 \langle  {  \bf{ l } }_k, \phi  \rangle = \sum_{i =1}^k   \int_{  \mathbb{R}^N} \big( V(|y|) - 1      \big)  \big( U_{ \overline{x}_i  } + U_{\underline{x}_i  }\big)  \phi +   \int_{  \mathbb{R}^N}  \Big(   \sum_{i =1}^k\big( U_{ \overline{x}_i  } + U_{\underline{x}_i  }\big)^p -  W_{r,h}^p    \Big) \phi.
\end{align}
Using the symmetric property for the functions,  we can know
\begin{align}\label{first}
& \sum_{i =1}^k   \int_{  \mathbb{R}^N} \big( V(|y|) - 1      \big)  \big( U_{ \overline{x}_i  } + U_{\underline{x}_i  }\big)  \phi   =2  k   \int_{  \mathbb{R}^N} \big( V(|y|) - 1      \big)   U_{ \overline{x}_1 }  \phi   \nonumber
\\[1mm]
& =  2  k   \Big(  \int_{  \mathbb{R}^N}  \big( V(|y+\overline{x}_1 |) - 1      \big)^2   U^2     \Big)^{\frac 12} \, \| \phi \|  \nonumber
\\[1mm]
& =  2  k   \Big(  \int_{  \mathbb{R}^N \setminus  B_{\delta_0 | \overline{x}_1 | } (0)     }  \big( V(|y+\overline{x}_1 |) - 1      \big)^2   U^2+           \int_{   B_{\delta_0 | \overline{x}_1 | } (0)     }  \big( V(|y+\overline{x}_1 |) - 1      \big)^2   U^2 \Big)^{\frac 12} \, \| \phi \|
 \nonumber
\\[1mm]
 &\leq k O_k( \frac{1}{r^m})  \|\phi\|
 \leq   \frac{C}{k^{  \frac{m-1}{2} + \sigma } }  \|\phi\|.
\end{align}
The last  inequality hold as $m>1$.

Note that explicitly
\begin{align*}
 U_{ \overline{x}_1  } + U_{\underline{x}_1 }  \geq   U_{ \overline{x}_i  } + U_{\underline{x}_i  } \quad \text{for } ~ y \in  \Omega_1,
\end{align*}
\begin{align*}
  U_{ \overline{x}_j  }    \geq   U_{\underline{x}_j  } \quad \text{for } ~ y \in  \Omega_1^+.
\end{align*}

For the second  term in \eqref{definlk},     we have
\begin{align} \label{second}
& \Big| \int_{  \mathbb{R}^N}  \Big(   \sum_{i =1}^k\big( U_{ \overline{x}_i  } + U_{\underline{x}_i  }\big)^p -  W_{r,h}^p    \Big) \phi \Big|  = k \Big| \int_{   \Omega_1}  \Big(   \sum_{i =1}^k\big( U_{ \overline{x}_i  } + U_{\underline{x}_i  }\big)^p -  W_{r,h}^p    \Big)   \phi  \Big|\nonumber
\\[1mm]
& \leq  k  \int_{   \Omega_1}  \big( U_{ \overline{x}_1  } + U_{\underline{x}_1 }\big)^{p-1}    \sum_{i =2}^k\big( U_{ \overline{x}_i  } + U_{\underline{x}_i  }\big) |\phi|   = 2 k  \int_{   \Omega_1^+ }  \big( U_{ \overline{x}_1  } + U_{\underline{x}_1 }\big)^{p-1}    \sum_{i =2}^k\big( U_{ \overline{x}_i  } + U_{\underline{x}_i  }\big)  |\phi|
\nonumber
\\[1mm]
& \leq  C k   \sum_{i =2}^k\int_{   \Omega_1^+ }   U_{ \overline{x}_1  }^{p-1}   U_{ \overline{x}_i  } |  \phi |    \leq       Ck  \sum_{i =2}^k    e^{ - \min{ \{  p-1-\tau,   1       \}        }    | \overline{x}_i - \overline{x}_1|   }    \Big(   \int_{   \Omega_1^+ }    |  \phi |^{p+1} \Big)^{\frac{1}{p+1}}
\nonumber
\\[1mm]
& \leq  C k^{ \frac{p}{p+1} }     \sum_{i =2}^k    e^{ - \min{ \{  p-1-\tau,   1       \}        }    | \overline{x}_i - \overline{x}_1|   }    \| \phi\|,
\end{align}
where $\tau$ is the constant small enough.

Noting  that for $(r,h)\in  \mathbb{S}_k$,  we can get
\begin{align*}
 \sum_{i =2}^k    e^{ - \min{ \{  p-1-\tau,   1       \}        }    | \overline{x}_i - \overline{x}_1|   }   & \leq C  e^{ - \min{ \{  p-1-\tau,   1       \}        }   \frac{2 \pi \sqrt{1-h^2 } }{k}  r } \nonumber
\\[1mm]
& \leq  \frac{C}{ k^{ \min{ \{  p-1-\tau,   1       \}        }   (m-\beta)  }      }.
\end{align*}
Combing  \eqref{second} and the following fact
\begin{align*}
\min{ \{   p-1, 1       \}        } m  - \frac{p}{p+1}> \frac{m-1}{2},
\end{align*}
we can obtain
\begin{align}\label{three}
\Big| \int_{  \mathbb{R}^N}  \Big(   \sum_{i =1}^k\big( U_{ \overline{x}_i  } + U_{\underline{x}_i  }\big)^p -  W_{r,h}^p    \Big) \phi \Big|  \leq  \frac{C}{k^{  \frac{m-1}{2} + \sigma } }  \|\phi\|.
\end{align}
  The result  of \eqref{estimatelk} follows from \eqref{first} and \eqref{three}. \qed
\medskip

In order to prove the Theorem \ref{main1}, we need the following proposition.  Readers can refer to \cite{LinNiWei},   \cite{wei&yan} for more details about the proposition.
\medskip
\begin{proposition} \label{pro2.4}
Suppose $ \Phi(r,h) = \phi_{r,h}(y)$ with $\Phi(r,h) $ be the map defined in   Proposition \ref{propos2}.  Define
\begin{align*}
 F(r,h) =   I(W_{r,h}  +     \phi_{r,h}(y) ),  \qquad  \forall  ~ (r,h)\in  \mathbb{S}_k.
 \end{align*}
If  $(r,h)$ is a critical point of $  F(r,h) $, then
 \begin{align*}
  u = W_{r,h}  +     \phi_{r,h}(y)
  \end{align*}
   is  a critical point
 $ I(u)$ in $ H^1(\mathbb{R}^N)$. \qed
\end{proposition}

%

%
{\textit{Proof of Proposition  \ref{pro2.4}}}:  The Proof  is similar to  the proof or Proposition \ref{func}. We omit it here.

\medskip
Now, we give the proof of Theorem \ref{main1}.

\medskip
{\textit{Proof of Theorem \ref{main1}  } }:     By Proposition \ref{pro2.4},  we need to show there  is $(r_k,h_k)\in \mathbb{S}_k$ which is a critical point of $F(r,h)$.

In fact,    from Proposition  \ref{func},   we can know
\begin{align*}
F(r,h) & =   I(W_{r,h} ) + {  \bf{ l } }  (\phi_{r,h}) + \frac 12  \langle{\bf{L}} \phi_{r,h},  \phi_{r,h}\rangle  -  {  \bf{ R } }(\phi_{r,h}) \nonumber
\\[1mm]
& =  I(W_{r,h} ) +  O_k(\|{\bf{l}}_k \|     \| \phi_{r,h}\| +      \|\phi_{r,h}\|^2  ) =  I(W_{r,h} )  +     O_k(\frac{1}{k^{m-1+\sigma}})\nonumber
\\[1mm]
& =
  k \Big(      \frac{A_1}{r^m} + A_2   - 2B_1     e^{-     2 \pi \sqrt{1-h^2 }    \frac{r}{k}    } -B_1 e^{-2    rh}   \nonumber
\\[1mm]
& \qquad + O_k( e^{-     2\pi (1+\sigma)\sqrt{1-h^2 }    \frac{r}{k}    }) +O_k(e^{-2 (1+\sigma)   rh}) +   O_k(\frac{1}{r^{m+\sigma}})  + O_k(\frac{1}{k^{m+\sigma}}) \Big),
\end{align*}
where $A_1,    A_2,        B_1$  are  constants  in \eqref{A1A2}.
Define  \[  F_1 ( r, h) =
  \frac{A_1}{r^m} + A_2   - 2B_1     e^{-     2 \pi \sqrt{1-h^2 }    \frac{r}{k}    } -B_1 e^{-2    rh}    \]

Then  we consider the system
 \begin{align}
 \begin{cases} \label{system1}
  F_{1,r}(r,h)  =     -   A_1    \frac{m }{r^{m+1 }}  + 4 B_1     \pi \sqrt{1-h^2 }  \frac{ e^{-     2 \pi \sqrt{1-h^2 }    \frac{r}{k}    } }{k}
      + 2  B_1    h e^{-2    rh}     =0,
 \\[5mm]
  F_{1,h}(r,h)  =     - 4B_1   \pi   \frac{hr}{ {\sqrt {1-h^2} }  }  \frac{ e^{-     2 \pi \sqrt{1-h^2 }    \frac{r}{k}    } }{k}   + 2B_1    r e^{-2    rh}  =0.
\end{cases}
\end{align}

From \eqref{system1}, we can get
\begin{align*}
 -   A_1    \frac{m }{r^{m+1 }} + 4B_1 \pi  \frac{ e^{-     2 \pi \sqrt{1-h^2 }    \frac{r}{k}    } }{k}  \big(\sqrt{1-h^2 }  +\frac{h^2}{ {\sqrt {1-h^2} }  }  \big)    =0.
\end{align*}
Define $$  H(r,h)= e^{-     2 \pi \sqrt{1-h^2 }    \frac{r}{k}    },     \quad  G(r,h)  = e^{-2    rh}.  $$
Then    \eqref{system1} implies

\begin{align}
 \begin{cases} \label{system2}
H(r,h)  =   \frac{      A_1    k \frac{m }{r^{m+1 }}      }{  4B_1 \pi      \big(\sqrt{1-h^2 }  +\frac{h^2}{ {\sqrt {1-h^2} }  }  \big)    }, \vspace{5mm} \\
G(r,h)  =  \frac{  2  \pi     h  }{    \sqrt {1-h^2} k  }  e^{-     2 \pi \sqrt{1-h^2 }    \frac{r}{k}    }.
\end{cases}
\end{align}
We define the space $\bar {\mathbb{S} }_k$:
\begin{align*}
\bar {\mathbb{S} }_k = \Big\{   ( {r}_k,   {h}_k)  \Big|  \quad
 {r}_k  = \big( \frac{m}{2\pi}+ o_k(1) \big)k \ln k, ~ \text{and }   {h}_k = \big( \frac{ \pi \, (m+2) }{  m  }   +o_k(1)  \big) \frac{1}{k} \Big\},
\end{align*}
and  the mapping
   $$ {\bf{T} } :   \bar {\mathbb{S} }_k  \rightarrow  \mathbb{R}^2   \quad \text{as} \quad  {\bf{T} }(r_k, h_k ) = \Big( H(r_k,h_k), G(r_k,h_k)  \Big).    $$
The system \eqref{system1} is equivalent to find a fixed point of
\begin{align}\label{system4}
(r, h ) &=    {\bf{T} } ^{-1}\Big(  \frac{        A_1  k \frac{m }{r^{m+1 }}      }{  4B_1 \pi   \big(\sqrt{1-h^2 }  +\frac{h^2}{ {\sqrt {1-h^2} }  }  \big)    },       \, \,     \frac{      h  }{  2 B_1  \sqrt {1-h^2}   } \,    \frac{          A_1    \frac{m }{r^{m+1 }}      }{     \big(\sqrt{1-h^2 }  +\frac{h^2}{ {\sqrt {1-h^2} }  }  \big)    }  \Big)\nonumber
\\[1mm] & = {\bf{A}}  (r,h) = :\Big( \,  {\bf{a}}_1  (r,h), \,  {\bf{a}}_2  (r,h)       \, \Big).
\end{align}
 in $ \bar {\mathbb{S} }_k $.
 By computing,  for $(r, h ) \in  \bar {\mathbb{S} }_k$, we have
 \begin{align}\label{system3}
 {\bf{A}}  (r,h)  &=  \big(1+ o_k(1) \big)   \Big( k\frac{(m+1) \ln {r} - \ln k }{2 \pi },         \, \, \,       \pi \, \frac{ \ln {h}- (m+1) \ln {r}  }{ k \big(   \ln k  - (m+1)   \ln {r}    \big) }      \Big) \in  \bar {\mathbb{S} }_k\nonumber
\\[1mm] & = \Big( \,  {\bf{a}}_1  (r,h), \,  {\bf{a}}_2  (r,h)       \, \Big).
 \end{align}
 And it is easy to show  that
 \begin{align*}
 | \,  {\bf{a}}_1   (r_1,h_1)  -  {\bf{a}}_1  (r_2,h_2)\,    | +  | \,  {\bf{a}}_2   (r_1,h_1)  -  {\bf{a}}_2  (r_2,h_2)\,    |  \leq \,  o_k(1)\,    \Big(  | r_1- r_2|+    \,  |h_1-h_2| \Big),
 \end{align*}
$\text{for all} ~  (r_1,h_1),   (r_2,h_2)\in \bar {\mathbb{S} }_k$.
By using the  Contraction Mapping principle we can prove that  there exist a  fixed point $( \bar{r}_k, \bar{h}_k )  \in \bar {\mathbb{S} }_k$ for   \eqref{system4}.  That's to say  $ F_1 ( r,h) $ have a critical point  $( \bar{r}_k, \bar{h}_k )  \in \bar {\mathbb{S} }_k$.

    Define
         \[    {   \bf M_1     } ( r, h )    = \left(
\begin{array}{cccc}
 F_{1,rr} & F_{1,rh}
\\[2mm]
F_{1,rh}& F_{1,hh}
\end{array}
\right). \]
By some simple calculations, we can know that  \[      F_{1,rr} \big| _{( r, h )   =( \bar{r}_k, \bar{h}_k )  }  < 0, \quad  F_{1,hh} \big| _{( r, h )   =( \bar{r}_k, \bar{h}_k )  }  < 0  \]
  and
  \[      F_{1,rr} \times   F_{1,hh} \big| _{( r, h )   =( \bar{r}_k, \bar{h}_k )  }  - F_{1,rh}^2 \big| _{( r, h )   =( \bar{r}_k, \bar{h}_k )  }   > 0.   \]
 So we can know that  $ ( \bar{r}_k, \bar{h}_k ) $ is a maximum point of    $ F_1 ( r,h) $.
  Then the
maximum of  $ F_1 ( r,h) $ in  $   {\mathbb{S} }_k$ can be achieved.
%
 Thus for the function $      F(r,h)$,   we can find a maximum  point $(r_k, h_k )   $ which is an interior point of $   {\mathbb{S} }_k$.  So $(r_k, h_k )   $ is a critical point of $      F(r,h)$.
Thus
$$  W_{r_k,h_k}  +     \phi_{r_k,h_k}(y)$$ is a critical point of $I(u)$.     This complete the proof of  Theorem \ref{main1}.
\qed

 \vspace{2mm}
\begin{remark} \label{remark6}
From \eqref{system1} and the assumptions $ h \in \big[ \alpha_0\frac{1}{k}, \alpha_1\frac{1}{k}\big]   $,  $ r\in [\beta_0 k \ln k, \beta_1 k \ln k]$, we can know that
\begin{align*}
e^{-2    rh} =    O_k \big(  \frac{h }{  k {\sqrt {1-h^2} }    }  \big)    e^{-     2 \pi \sqrt{1-h^2 }    \frac{r}{k}    } = O_k(\frac{1}{k^2})  e^{-     2 \pi \sqrt{1-h^2 }    \frac{r}{k}    }.
\end{align*}
For the  solvability  system \eqref{system1}, we must have
\begin{align*}
e^{-     2 \pi \sqrt{1-h^2 }    \frac{r}{k}    }= O_k(\frac{1}{k^{m}})  \quad \text{and} \quad   e^{-2    rh}  =  O_k(\frac{1}{k^{m+2}}).
\end{align*}
That's the reason why we make the assumption \eqref{H2} for $(r,h)$.
\end{remark}

 \vspace{3mm}
%
\section{ Some useful estimates and Lemmas}\label{appendixA}
We first  give  some useful estimates and Lemmas.

\begin{lemma}\label{lemma1}
For  $r,h $  be   the parameters   in \eqref{overunder}  and  any $\eta \in (0,1] $, there  is a constant $ C>0$, such that
\begin{align*}
\sum_{i =2}^k  U_{ \overline{x}_j} (y)& \leq C e^{- \eta \sqrt{1-h^2} r\frac{\pi}{k}}e^{   -(1-\eta) |y-\overline{x}_1| },  \quad \text{for all} ~  y \in \Omega_1^+,
\end{align*}
and
\begin{align*}
\sum_{i =2}^k  U_{ \underline{x}_j} (y)& \leq C e^{- \eta \sqrt{1-h^2} r\frac{\pi}{k}} e^{   -(1-\eta) |y-\overline{x}_1| },  \quad \text{for all} ~  y \in \Omega_1^+,
\end{align*}
\begin{align*}
 U_{ \underline{x}_1} (y)\leq C   e^{ -   \eta  hr}   e^{-(1-\eta) | y -   \overline{x}_1  | },   \quad \text{for all} ~  y \in \Omega_1^+.
\end{align*}
\end{lemma}
{\textit{Proof of Lemma \ref{lemma1}  } }:     The proof  of Lemma \ref{lemma1} is similar to  Lemma A.1 in \cite{wei&yan}.
\qed
\medskip

In this appendices, we assume  $(r,h)  \in \mathbb{S}_k $, where    $\mathbb{S}_k$ is defined in \eqref{H2}.

%

\medskip
\begin{lemma}  \label{lemma2}
For $ i = 2, \cdots, k, $  there exist some small constant $\sigma>0$ such that
 the following expansions  hold
 \begin{align} \label{A1}
 \int_{  \mathbb{R}^N}  U^p_{     \overline{x}_1  }   U_{\overline{x}_i  }    =      B_1e^{-|\overline{x}_1 -      \overline{x}_i| }   +   O_k( e^{-(1+\sigma) | \overline{x}_1 -      \overline{x}_i| }   ),
 \end{align}
 \begin{align}\label{A2}
  \int_{  \mathbb{R}^N}  U^p_{\underline{x}_i }  U_{\overline{x}_1  }   =   B_1e^{-|\overline{x}_1 -      \underline{x}_i| }   +   O_k( e^{-(1+\sigma) | \overline{x}_1 -      \underline{x}_i| }   ),
 \end{align}
 \begin{align}\label{A3}
 \int_{  \mathbb{R}^N}     U^p_{\underline{x}_1 }  U_{\overline{x}_1  }   =  B_1e^{-2rh} + O_k( e^{-2(1+\sigma) rh} ),
 \end{align}
 where $B_1$ is  defined in \eqref{A1A2}.
 \end{lemma}
{\textit{Proof of Lemma \ref{lemma2}  } }:       By calculation, we have
 \begin{align} \label{sum1}
        \int_{  \mathbb{R}^N}  U^p_{     \overline{x}_1  }   U_{\overline{x}_i  }  & = \int_{  \mathbb{R}^N}  U^p(z)    U(z+    \overline{x}_1 -      \overline{x}_i         ) \nonumber
   \\[1mm]
&
 \overset{\text{     $w_i = \overline{x}_1 -      \overline{x}_i $     }}{=}  \int_{  \mathbb{R}^N}  U^p(z)    U(z+    w_i        )\nonumber
   \\[1mm]
   & = \int_{  \mathbb{R}^N   \setminus    B_{ \delta |w_i|}  (0) }  U^p(z)    U(z+    w_i        )  +    \int_{    B_{ \delta |w_i| }(0) }  U^p(z)    U(z+    w_i        ).
\end{align}
with $\delta>0$  is small constant.  Then we have
\begin{align}\label{pro1}
 \int_{  \mathbb{R}^N   \setminus    B_{ \delta |w_i|}  (0) }  U^p(z)    U(z+    w_i        ) &\leq C e^{-|w_i|}  \int_{  \mathbb{R}^N   \setminus    B_{ \delta |w_i|}  (0) }     e^{- ( p- 1) |z| }  \nonumber
   \\[1mm]
&    =  O_k( e^{-(1+\sigma) | w_i| }   ) =  O_k( e^{-(1+\sigma) | \overline{x}_1 -      \overline{x}_i| }   ).
\end{align}
We have the expansion in $  B_{ \delta |w_i| }(0)$
 $$  \big| z+   | w_i | e_1\big|= |w_i| \big(   1+  \frac{1}{|w_i|  }  e_1z +     O_k(   \frac{|z|^2}{|w_i| ^2}   ) \big),     $$
 where $e_1= (1, \cdots, 0) $,
 then  using the symmetry of $U$,
\begin{align}\label{pro2}
  \int_{    B_{ \delta |w_i| }(0) }  U^p(z)    U(z+    w_i       )     &=\int_{    B_{ \delta |w_i| }(0) }  U^p(z)    U(z+    |w_i|e_1        )
 =    \int_{    B_{ \delta |w_i| }(0) }  U^p(z)    e^{-   \big| z+   | w_i | e_1\big|  }  \nonumber
    \\[1mm]
&  =     \int_{    B_{ \delta |w_i| }(0) }  U^p(z)    e^{-|w_i| \big(   1+  \frac{1}{|w_i|  }  e_1z +     O_k(   \frac{|z|^2}{|w_i| ^2}   ) \big) }   \nonumber
    \\[1mm]
&  =  e^{-|\overline{x}_1 -      \overline{x}_i| }   \int_{    \mathbb{R}^N  }  U^p(z)    e^{-         z_1    }     (1+     o_k(1)  ).
\end{align}
Now combining \eqref{sum1}, \eqref{pro1} and \eqref{pro2}, we can get
 \begin{align*}
 \int_{  \mathbb{R}^N}  U^p_{     \overline{x}_1  }   U_{\overline{x}_i  }    =      B_1e^{-|\overline{x}_1 -      \overline{x}_i| }   +   O_k( e^{-(1+\sigma) | \overline{x}_1 -      \overline{x}_i| }   ),
 \end{align*}
 with $B_1= \int_{    \mathbb{R}^N  }  U^p(z)    e^{-         z_1    }    $.    Similarity, we can  easily obtain \eqref{A2}, \eqref{A3}.
Thus we complete the proof of Lemma \ref{lemma2}.

\qed

\begin{lemma}\label{lemma3}
There exist some small constant $\sigma>0$   such that the following expansions hold
\begin{align}  \label{sum2}
\sum_{i =2}^k  \int_{  \mathbb{R}^N}  U^p_{     \overline{x}_1  }   U_{\overline{x}_i  }    =     2 B_1     e^{  - 2\pi  \sqrt{1-h^2}  \frac r k                }     +    O_k( e^{  - 2(1+\sigma) \pi  \sqrt{1-h^2}  \frac r k                }  ),
\end{align}
\begin{align}\label{sum3}
   \sum_{i =1}^k \int_{  \mathbb{R}^N}    U^p_{\underline{x}_j }  U_{\overline{x}_1 } & =  B_1   e^{-2    rh}        +     O_k( e^{- 2(1+\sigma)    rh}  )  +   O_k(     e^{-     2 \pi (1+\sigma )\sqrt{1-h^2 }    \frac{r}{k}    }   ).
\end{align}
\end{lemma}
{\textit{Proof of Lemma \ref{lemma3}  } }:
Recalling the definitions $ \overline{x}_j, \underline{x}_j$, we can know
\begin{align*}
   |  \overline{x}_j- \overline{x}_1|^2  = 4r^2  ( 1- h^2)    \sin^2 { \frac{ ( j-1 )\pi  }{k}    },
\end{align*}
 \begin{align*}
   |  \overline{x}_1- \underline{x}_1|^2=4  r^2 h^2,
 \end{align*}
 and
\begin{align*}
 |  \underline{x}_j- \overline{x}_1|^2  =4  r^2 \Big[  ( 1- h^2)    \sin^2 { \frac{ ( j-1 )\pi  }{k}    }   +h^2     \Big].
\end{align*}
From  Lemma \ref{lemma2},  we can get
 \begin{align} \label{pro3}
   k \sum_{i =2}^k  \int_{  \mathbb{R}^N}  U^p_{     \overline{x}_1  }   U_{\overline{x}_i  }  & =   k  B_1 \sum_{i =2}^k  e^{-|\overline{x}_1 -      \overline{x}_i| }   +  k \sum_{i =2}^k  O_k( e^{-(1+\sigma) | \overline{x}_1 -      \overline{x}_i| }   ),
 \end{align}
 with $B_1= \int_{    \mathbb{R}^N  }  U^p(z)    e^{-         z_1    }    $.
 Without  lose  of generality, we can assume that  $k$ is even.
 Then
  \begin{align} \label{pro6}
    B_1 \sum_{j =2}^k  e^{-|\overline{x}_1 -      \overline{x}_i| }
&=      B_1   \sum_{i =3}^{ k/2}   e^{  - 2r    \sqrt{1-h^2}      \sin{\frac{( j-1 )\pi } {k }            }       }
\nonumber
\\[1mm] & \quad
+ B_1   \sum_{i = { k/2}  + 1}^{k-1}  e^{  - 2r    \sqrt{1-h^2}      \sin{\frac{( j-1 )\pi } {k }            }       }
+  2  B_1e^{  - 2r    \sqrt{1-h^2}    \sin{\frac{ \pi } {k }            }         }.
 \end{align}
 Consider   \begin{align*}
 c_3 \frac{ ( j-1 )\pi  }{k}      \leq   \sin { \frac{ ( j-1 )\pi  }{k}    }     \leq c_4 \frac{ ( j-1 )\pi  }{k}, \quad \text{for} ~j \in           \big\{ 3, \cdots,    \frac k2 \big\},
 \end{align*}
 with $ \frac 12  <  c_3 \leq c_4\leq 1.$
Then we have
\begin{align}\label{pro4}
  \sum_{j=3}^{ k/2}   e^{  - 2r    \sqrt{1-h^2}      \sin{\frac{( j-1 )\pi } {k }            }       }
& \leq    \sum_{j=3}^{ k/2}   e^{  - 2r    \sqrt{1-h^2}     \frac{ c_3( j-1 )\pi } {k }                 }      \nonumber
    \\[1mm]  & =   \frac{   e^{  - 4r    \sqrt{1-h^2}     \frac{ c_3\pi } {k }                     }   -      e^{  - 2r    \sqrt{1-h^2}     \frac{ c_3  (\frac k2 -3) \pi } {k }                     }          }{       (    1-  e^{  - 2r    \sqrt{1-h^2}     \frac{ c_3 \pi } {k }                     }   )     }    \nonumber
   \\[1mm]
&  = O_k( e^{  - 2(1+\sigma) \pi  \sqrt{1-h^2}  \frac r k                }  ).
\end{align}
Using symmetry  of  function $\sin x$,  we can easily show
  $$ B_1   \sum_{j = { k/2}  + 1}^{k-1}  e^{  - 2r    \sqrt{1-h^2}      \sin{\frac{( j-1 )\pi } {k }            }       } = O_k( e^{  - 2(1+\sigma) \pi  \sqrt{1-h^2}  \frac r k                }  ),  $$ in a same  manner as  \eqref{pro4}.
Next,  we have
\begin{align}\label{pro5}
 e^{  - 2r    \sqrt{1-h^2}    \sin{\frac{ \pi } {k }            }         }  &=  e^{  - 2r    \sqrt{1-h^2}   (  \frac{ \pi } {k }      + O_k( \frac{ \pi^3 } {k^3 }  )                } \nonumber
   \\[1mm]
&  = e^{  - 2\pi  \sqrt{1-h^2}  \frac r k                }e^{  - 2r    \sqrt{1-h^2}       O_k( \frac{ \pi^3 } {k^3 }  )                }\nonumber
   \\[1mm]
&  = e^{  - 2\pi  \sqrt{1-h^2}  \frac r k                }   +  O_k( e^{  - 2(1+\sigma) \pi  \sqrt{1-h^2}  \frac r k                }  ).
\end{align}
Combining \eqref{pro3}, \eqref{pro6}, \eqref{pro4} and \eqref{pro5},   we can get
\begin{align}
   \sum_{i =2}^k  \int_{  \mathbb{R}^N}  U^p_{     \overline{x}_1  }   U_{\overline{x}_i  }    =     2   B_1     e^{  - 2\pi  \sqrt{1-h^2}  \frac r k                }     +    O_k( e^{  - 2(1+\sigma) \pi  \sqrt{1-h^2}  \frac r k                }  ).
 \end{align}

 Next, from \eqref{A2} and \eqref{A3}, we have
\begin{align}\label{mathbbI13}
&     \sum_{j =1}^k \int_{  \mathbb{R}^N}    U^p_{\underline{x}_j }  U_{\overline{x}_1 } \nonumber
   \\[1mm]
&  =   B_1     e^{-     |  \overline{x}_1- \underline{x}_1| }   +   B_1 \sum_{j =2}^k  e^{-     |  \underline{x}_j- \overline{x}_1| }  +   O_k( e^{- (1+\sigma)    |  \overline{x}_1- \underline{x}_1| }  )+   O    (\sum_{j =2}^k   e^{-   (1+\sigma)   |  \underline{x}_j- \overline{x}_1| }     )\nonumber
   \\[1mm]
&  =  B_1   e^{-2    rh} +   B_1\sum_{j =2}^k  e^{-    2  r \Big[  ( 1- h^2)    \sin^2 { \frac{ ( j-1 )\pi  }{k}    }   +h^2     \Big]^{\frac 12} } \nonumber
   \\[1mm]
& \quad +  O_k( e^{- 2(1+\sigma)    rh}  ) +   O \Big( \sum_{j =2}^k  e^{-    2(1+\sigma)  r \Big[  ( 1- h^2)    \sin^2 { \frac{ ( j-1 )\pi  }{k}    }   +h^2     \Big]^{\frac 12} }\Big) \nonumber
   \\[1mm]
&  =     B_1   e^{-2    rh} + 2  B_1  e^{-    2  r \Big[  ( 1- h^2)    \sin^2 { \frac{   \pi  }{k}    }   +h^2     \Big]^{\frac 12} }  +  O_k( e^{- 2(1+\sigma)    rh}  ) + O_k( e^{  - 2(1+\sigma) \pi  \sqrt{1-h^2}  \frac r k                }  ).
\end{align}
Recalling  $(r,h)  \in \mathbb{S}_k$,      we have
\begin{align}\label{pro7}
e^{-    2  r \Big[  ( 1- h^2)    \sin^2 { \frac{   \pi  }{k}    }   +h^2     \Big]^{\frac 12} }    & =     e^{-    2  r        ( 1- h^2)^{\frac 12}   \sin { \frac{ \pi  }{k}    } \Big[  1  +   \frac{h^2 }{  ( 1- h^2)    \sin^2 { \frac{  \pi  }{k}    } }        \Big]^{\frac 12} }       \nonumber
   \\[1mm]
&=  O_k( e^{-    2 (1+\sigma )  r        ( 1- h^2)^{\frac 12}     { \frac{ \pi  }{k}    } }).
\end{align}

Then together with \eqref{mathbbI13} and \eqref{pro7},   we can know easily
\begin{align}
   \sum_{j =1}^k \int_{  \mathbb{R}^N}    U^p_{\underline{x}_j }  U_{\overline{x}_1 } & =  B_1   e^{-2    rh}        +     O_k( e^{- 2(1+\sigma)    rh}  )  +   O_k(     e^{-     2 \pi (1+\sigma )\sqrt{1-h^2 }    \frac{r}{k}    }   ).
\end{align}
We complete the proof of Lemma \ref{lemma3}.
\qed

\bigskip
\section{ Proof of Proposition \ref{func} }\label{appendixB}

In this section, we will turn to  the proof of Proposition \ref{func}.

{\textit{Proof of Proposition \ref{func}  } }:
Now we calculate $ I(W_{r,h} )$:

\begin{align}\label{energyexpan}
 I(W_{r,h} ) &= \frac 12 \int_{  \mathbb{R}^N}   \Big\{  |  \nabla  W_{r,h} |^2 + V(|y|) W_{r,h}^2  \Big\} - \frac{1}{p+1} \int_{\mathbb{R}^N}  |W_{r,h}|^{p+1} \nonumber
 \\[1mm]
&  = \frac 12 \int_{  \mathbb{R}^N}   \Big\{ |  \nabla  W_{r,h} |^2  +   | W_{r,h}| ^2 \Big\}  \nonumber
\\[1mm]
   &\quad + \frac 12 \int_{  \mathbb{R}^N}  \big( V(|y|)- 1 \big )  | W_{r,h}| ^2 - \frac{1}{p+1} \int_{\mathbb{R}^N}  |W_{r,h}|^{p+1}  \nonumber
   \\[1mm]
  &  = \mathbb{I}_1 +  \mathbb{I}_2 + \mathbb{I}_3.
\end{align}

By using the symmetry and the equation for $U$, we can get
 \begin{align}\label{mathbbI1}
 \mathbb{I}_1 &= \frac 12 \int_{  \mathbb{R}^N}   \big( \{ |  \nabla  W_{r,h} |^2  +   | W_{r,h}| ^2 \big)    = \frac 12 \int_{  \mathbb{R}^N}    \big (   - \Delta W_{r,h} + W_{r,h}  \big ) W_{r,h} \nonumber
 \\[1mm]
&  = \frac 12 \int_{  \mathbb{R}^N} \sum_{j=1}^k   \big (    - \Delta  U_{\overline{x}_j  } + U_{\overline{x}_j  }       - \Delta U_{\underline{x}_j  }  + U_{\underline{x}_j  }   \big) \cdot  \sum_{i =1}^k  \big ( U_{\overline{x}_i } + U_{\underline{x}_i }    \big) \nonumber
 \\[1mm]
&  = \frac 12   \sum_{j=1}^k  \sum_{ j =1}^k   \int_{  \mathbb{R}^N}     \big (   U^p_{\overline{x}_j  }
  +   U^p_{\underline{x}_j  }     \big) \cdot    \big ( U_{\overline{x}_i } + U_{\underline{x}_i }    \big)  \nonumber
   \\[1mm]
&  = \frac 12   \sum_{j=1}^k  \sum_{j =1}^k   \int_{  \mathbb{R}^N}     \big (    U^p_{\overline{x}_j  }   U_{\overline{x}_i  }
  +   U^p_{\underline{x}_j }  U_{\overline{x}_i  }  + U^p_{\overline{x}_j  }  U_ { \underline{x}_i  } + U^p_{\underline{x}_j }  U_{\underline{x}_i  }      \big)  \nonumber
   \\[1mm]
&  =   k  \int_{  \mathbb{R}^N}        U^{p+1} + k \sum_{i =2}^k  \int_{  \mathbb{R}^N}  U^p_{\overline{x}_1  }   U_{\overline{x}_i  }  +  k   \sum_{j=1}^k \int_{  \mathbb{R}^N}    U^p_{\underline{x}_j }  U_{\overline{x}_1  }.
 \end{align}


Next we calculate $\mathbb{I}_2$, using symmetry and  Lemma \ref{lemma1}
\begin{align}\label{mathbbI2}
\mathbb{I}_2 & = \frac 12 \int_{  \mathbb{R}^N}  \big( V(|y|)- 1 \big )  | W_{r,h}| ^2
 \nonumber
   \\[1mm]
&  =      k   \int_{  \Omega_1^+}    \big( V(|y|)- 1 \big )   \big(    U_{ \underline{x}_1} + U_{ \overline{x}_1}   + \sum_{j =2}^k  U_{ \underline{x}_j} +\sum_{j =2}^k  U_{ \overline{x}_j}  \big)^2 \nonumber
   \\[1mm]
&  =     k   \int_{  \Omega_1^+}    \big( V(|y|)- 1 \big )   \Big( U_{ \overline{x}_1} + O_k\big(  e^{ -\frac 12 hr}   e^{-\frac 12 | y -   \overline{x}_1  | }    \nonumber
   \\[1mm]  &\quad
+   e^{- \frac 12 \sqrt{1-h^2} r\frac{\pi}{k}}e^{   -\frac 12 |y-\overline{x}_1| }+       e^{- \frac 12\sqrt{1-h^2} r\frac{\pi}{k}}  e^{   -\frac 12 |y-\overline{x}_1| }   \big) \Big)^2\nonumber
   \\[1mm]
&  =    k   \int_{  \Omega_1^+}  \big( V(|y|)- 1 \big )   U^2_{ \overline{x}_1}
 + k   O_k\Big(   \int_{  \Omega_1^+}  \big( V(|y|)- 1 \big )            e^{   -\frac 12 |y-\overline{x}_1| }      U_{ \overline{x}_1} \Big)  \nonumber
   \\[1mm]
&  =  k \Big\{   \frac{A_1}{r^m}+ O_k(\frac{1}{r^{m+\sigma}})\Big\},
\end{align}
where $A_1   $ is defined in \eqref{A1A2}  and the last  equality holds duo to the asymptotic expression of $V(y)$.

We consider $ \mathbb{I}_3 $.   Suppose $ p\leq 3 $, then for $ y\in \Omega_1^+$,
\begin{align*}
 |W_{r,h}|^{p+1} & = \Big( U_{\underline{x}_1 }+ U_{\overline{x}_1} +  \sum_{j=2}^kU_{\overline{x}_j }+ \sum_{j=2}^k U_{\underline{x}_j }  \Big) ^{p+1} \nonumber
   \\[1mm]
&  =    U_{\overline{x}_1}^{p+1}+ (p+1)   U_{\overline{x}_1}^{p}   \big( U_{\underline{x}_1 }+ \sum_{j=2}^kU_{\overline{x}_j }+ \sum_{j=2}^k U_{\underline{x}_j } \big) \nonumber
   \\[1mm]
    & \quad +  O_k\Big(  U_{\overline{x}_1}^{\frac{p+1}{2}}   \big( U_{\underline{x}_1 }+ \sum_{j=2}^kU_{\overline{x}_j }+ \sum_{j=2}^k U_{\underline{x}_j } \big)^{\frac{p+1}{2}} \Big).
 \end{align*}

   Using the Lemma \eqref{lemma1},   we can get  for $ y\in \Omega_1^+$,
 \begin{align*}
 &U_{\overline{x}_1}^{\frac{p+1}{2}}   \big( U_{\underline{x}_1 }+ \sum_{j=2}^kU_{\overline{x}_j }+ \sum_{j=2}^k U_{\underline{x}_j } \big)^{\frac{p+1}{2}}  \nonumber
   \\[1mm]
   & = U_{\overline{x}_1}^{\frac{p+1}{2}}          \big( U_{\underline{x}_1 }+ \sum_{j=2}^kU_{\overline{x}_j }+ \sum_{j=2}^k U_{\underline{x}_j } \big)  \big( U_{\underline{x}_1 }+ \sum_{j=2}^kU_{\overline{x}_j }+ \sum_{j=2}^k U_{\underline{x}_j } \big)^{\frac{p-1}{2}} \nonumber
   \\[1mm]
   & \leq   C  U_{\overline{x}_1}^{(1-\eta) p }    \Big(     e^{ - hr}   e^{-(1-\eta) | y -   \overline{x}_1  | } +   e^{- \eta \sqrt{1-h^2} r\frac{\pi}{k}}  e^{   -(1-\eta) |y-\overline{x}_1| }   \Big)^{\frac{p-1}{2}}  \nonumber
   \\[1mm]
   &     \qquad    \cdot  \big( U_{\underline{x}_1 }+ \sum_{j=2}^kU_{\overline{x}_j }+ \sum_{j=2}^k U_{\underline{x}_j } \big).
 \end{align*}
  For  any $ (r,h)\in S_k, $ and $ y\in \Omega_1^+$  we have
\begin{align*}
U_{\overline{x}_1  }    U_{\underline{x}_1 } \leq &C e^{- | \overline{x}_1  -  \underline{x}_1|      }   = C e^{-2hr },
\end{align*}
\begin{align*}
  U_{\overline{x}_1 }  \sum_{j=2}^k U_{\overline{x}_j  }   \leq C \sum_{j=2}^k e^{- | \overline{x}_1  -  \overline{x}_j|      }   \leq C e^{-2\pi  \sqrt{1-h^2}  \frac r k }+ Ce^{-2(1+\sigma) \pi  \sqrt{1-h^2}  \frac r k },
\end{align*}
and
\begin{align*}
U_{\overline{x}_1  }     \sum_{j=2}^k U_{\underline{x}_j } \leq C\sum_{j=2}^k  e^{- | \overline{x}_1  -  \underline{x}_j|      }   \leq C e^{-2\pi  \sqrt{1-h^2}  \frac r k }.
\end{align*}
 So, we obtain that for $ p\in (1,3], $
 \begin{align}  \label{mathbbI3}
 \mathbb{I}_3& = - \frac{1}{p+1} \int_{\mathbb{R}^N}  |W_{r,h}|^{p+1}  \nonumber
   \\[1mm]  & = - \frac{2k}{p+1}  \int_{ \Omega_1^+}  \Big( U_{\underline{x}_1 }+ U_{\overline{x}_1} +  \sum_{j=2}^kU_{\overline{x}_j }+ \sum_{j=2}^k U_{\underline{x}_j }  \Big) ^{p+1} \nonumber
   \\[1mm]  & = - \frac{2k}{p+1}   \int_{ \Omega_1^+}    \Big(  U_{\overline{x}_1}^{p+1}+ (p+1)   U_{\overline{x}_1}^{p}   \big( U_{\underline{x}_1 }+ \sum_{j=2}^kU_{\overline{x}_j }+ \sum_{j=2}^k U_{\underline{x}_j } \big)   \Big) \nonumber
   \\[1mm]
    & \quad  - \frac{2k}{p+1}        \int_{ \Omega_1^+}   O_k\Big(  U_{\overline{x}_1}^{\frac{p+1}{2}}   \big( U_{\underline{x}_1 }+ \sum_{j=2}^kU_{\overline{x}_j }+ \sum_{j=2}^k U_{\underline{x}_j } \big)^{\frac{p+1}{2}} \Big)   \nonumber
   \\[1mm]  & =   - \frac{2k}{p+1}   \int_{ \mathbb{R}^N }    \Big(  U_{\overline{x}_1}^{p+1}+ (p+1)   U_{\overline{x}_1}^{p}   \big( U_{\underline{x}_1 }+ \sum_{j=2}^kU_{\overline{x}_j }+ \sum_{j=2}^k U_{\underline{x}_j } \big)   \Big)  \nonumber
   \\[1mm]  & \quad
   +   k \Big(     e^{ - hr}    +   e^{- \eta \sqrt{1-h^2} r\frac{\pi}{k}}  +  e^{- \eta \sqrt{1-h^2} r\frac{\pi}{k}} \Big)^{\frac{p-1}{2}} \nonumber
   \\[1mm]  & \quad
    \cdot O_k( e^{- | \overline{x}_1  -  \underline{x}_1|      }  + \sum_{j=2}^k e^{- | \overline{x}_1  -  \overline{x}_j|      }     +   \sum_{j=2}^k e^{- | \overline{x}_1  -  \underline{x}_j|   })
    \nonumber
   \\[1mm]  & =   - \frac{2k}{p+1}   \int_{ \mathbb{R}^N }    \Big(  U_{\overline{x}_1}^{p+1}+ (p+1)   U_{\overline{x}_1}^{p}   \big( U_{\underline{x}_1 }+ \sum_{j=2}^kU_{\overline{x}_j }+ \sum_{j=2}^k U_{\underline{x}_j } \big)   \Big)
   \nonumber
   \\[1mm]  & \quad  +   k O_k( e^{-2(1+\sigma)\pi  \sqrt{1-h^2}  \frac r k }  + e^{-2(1+\sigma) r h }   ).
   \end{align}

Now we consider $ p>3$. Then for any $ y\in \Omega_1^+$,
\begin{align*}
 |W_{r,h}|^{p+1} & =\Big( U_{\underline{x}_1 }+ U_{\overline{x}_1} +  \sum_{j=2}^kU_{\overline{x}_j }+ \sum_{j=2}^k U_{\underline{x}_j }  \Big) ^{p+1}
 \nonumber
   \\[1mm]
&  =   U_{\overline{x}_1}^{p+1}+ (p+1)   U_{\overline{x}_1}^{p}   \big( U_{\underline{x}_1 }+ \sum_{j=2}^kU_{\overline{x}_j }+ \sum_{j=2}^k U_{\underline{x}_j } \big) \nonumber
   \\[1mm]
    & \quad +  O_k\Big(  U_{\overline{x}_1}^{p-1}   \big( U_{\underline{x}_1 }+ \sum_{j=2}^kU_{\overline{x}_j }+ \sum_{j=2}^k U_{\underline{x}_j } \big)^{2} \Big).
\end{align*}
Since $p-1>2$,  similar  to the proof of \eqref{mathbbI3},  we have
\begin{align}\label{mathbbI31}
\mathbb{I}_3& = - \frac{1}{p+1} \int_{\mathbb{R}^N}  |W_{r,h}|^{p+1}  \nonumber
   \\[1mm]
    & =
 - \frac{2k}{p+1}  \int_{ \Omega_1^+}  \Big( U_{\underline{x}_1 }+ U_{\overline{x}_1} +  \sum_{j=2}^kU_{\overline{x}_j }+ \sum_{j=2}^k U_{\underline{x}_j }  \Big) ^{p+1} \nonumber
   \\[1mm]
    & =   - \frac{2k}{p+1}   \int_{ \mathbb{R}^N }    \Big(  U_{\overline{x}_1}^{p+1}+ (p+1)   U_{\overline{x}_1}^{p}   \big( U_{\underline{x}_1 }+ \sum_{j=2}^kU_{\overline{x}_j }+ \sum_{j=2}^k U_{\underline{x}_j } \big)   \Big)  \nonumber
   \\[1mm]
    & \quad
    +
    kO_k\Big(     \int_{ \Omega_1^+}    U_{\overline{x}_1}^{p-1}   \big( U_{\underline{x}_1 }+ \sum_{j=2}^kU_{\overline{x}_j }+ \sum_{j=2}^k U_{\underline{x}_j } \big)^{2} \Big) \nonumber
   \\[1mm]
    & =
      - \frac{2k}{p+1}   \int_{ \mathbb{R}^N }    \Big(  U_{\overline{x}_1}^{p+1}+ (p+1)   U_{\overline{x}_1}^{p}   \big( U_{\underline{x}_1 }+ \sum_{j=2}^kU_{\overline{x}_j }+ \sum_{j=2}^k U_{\underline{x}_j } \big)   \Big)  \nonumber
   \\[1mm]  & \quad
   +   k \Big(     e^{ - hr}    +   e^{- \eta \sqrt{1-h^2} r\frac{\pi}{k}}   +  e^{- \eta \sqrt{1-h^2} r\frac{\pi}{k}} \Big)
    \cdot O_k( e^{-2hr } + e^{-2\pi  \sqrt{1-h^2}  \frac r k }    ) \nonumber
   \\[1mm]
    & =
      - \frac{2k}{p+1}   \int_{ \mathbb{R}^N }    \Big(  U_{\overline{x}_1}^{p+1}+ (p+1)   U_{\overline{x}_1}^{p}   \big( U_{\underline{x}_1 }+ \sum_{j=2}^kU_{\overline{x}_j }+ \sum_{j=2}^k U_{\underline{x}_j } \big)   \Big)  \nonumber
   \\[1mm]
    & \quad
      +
  k O_k( e^{-2(1+\sigma) \pi \sqrt{1-h^2}  \frac  r k }  + e^{-2(1+\sigma) r h }   ).
    \end{align}
Combining \eqref{energyexpan},  \eqref{mathbbI1}, \eqref{mathbbI2}, \eqref{mathbbI3}, \eqref{mathbbI31}, we can get
\begin{align}\label{Wrh}
 I(W_{r,h} )
  &  = \mathbb{I}_1 +  \mathbb{I}_2 + \mathbb{I}_3  \nonumber
   \\[1mm]
    &  =
    k  \int_{  \mathbb{R}^N}        U^{p+1} + k \sum_{i =2}^k  \int_{  \mathbb{R}^N}  U^p_{\overline{x}_1  }   U_{\overline{x}_i  } + k \sum_{j =1}^k\int_{  \mathbb{R}^N}  U^p_{\underline{x}_j }  U_{\overline{x}_1  }   \nonumber
   \\[1mm]
    & \quad
      - \frac{2k}{p+1}   \int_{ \mathbb{R}^N }    \Big(  U^{p+1}+ (p+1)   U_{\overline{x}_1}^{p}   \big( U_{\underline{x}_1 }+ \sum_{j=2}^kU_{\overline{x}_j }+ \sum_{j=2}^k U_{\underline{x}_j } \big)   \Big)  +    k \Big\{   \frac{A_1}{r^m}+ O_k( \frac{1}{r^{m+\sigma}})\Big\}\nonumber
   \\[1mm]
    &  = k \Big\{   \frac{A_1}{r^m}+ A_2+  O_k( \frac{1}{r^{m+\sigma}})\Big\}-  k \sum_{i =2}^k  \int_{   \mathbb{R}^N} U^p_{\overline{x}_1  }   U_{\overline{x}_i  }    -k \sum_{j =1}^k\int_{  \mathbb{R}^N}  U^p_{\underline{x}_j }  U_{\overline{x}_1  }\nonumber
   \\[1mm]
    &  \quad
    +
  k O_k( e^{-2(1+\sigma)\pi  \sqrt{1-h^2}  \frac r k }  + e^{-2(1+\sigma) r h }   ).
\end{align}
Combining    Lemma \ref{lemma3} and \eqref{Wrh},   we have
\begin{align}
 I(W_{r,h} )
&   =   k \Big(      \frac{A_1}{r^m} + A_2    - 2 B_1    e^{-     2 \pi \sqrt{1-h^2 }    \frac{r}{k}    }  -  {B}_1  k  e^{-2    rh} \nonumber
   \\[1mm]   & \qquad    + O_k( \frac{1}{r^{m+\sigma}})  + O_k(     e^{- 2 \pi (1+\sigma)      \sqrt{1-h^2 }      \frac{r}{k}             }        )  + O_k( e^{- 2(1+\sigma)    rh}  ) \Big),
\end{align}
 where    $A_1, A_2, B_1$  are  defined in \eqref{A1A2}. \qed

\bigskip
 \vspace{3mm}

{\bf Acknowledgements: }
 L. Duan   was supported by the China Scholarship Council and NSFC grants (No.11771167). M. Musso has been supported by EPSRC research
Grant EP/T008458/1.
This paper was completed during the visit L. Duan   to  M. Musso  at the  University of Bath.  L. Duan  would like to thank the Department of Mathematical Sciences  for its warm hospitality and supports.

 \end{document}